\documentclass[12pt]{article} 
\usepackage[english]{babel}
\usepackage[latin1]{inputenc} 
\usepackage[T1]{fontenc}
\usepackage{color,graphicx} 
\usepackage{amssymb} 
\usepackage{amsmath}
\usepackage{amsfonts,txfonts,pxfonts,array} 
\usepackage{bbm}
\usepackage{hyperref} 
\usepackage[nottoc, notlof, notlot]{tocbibind}
\usepackage{setspace} 
\usepackage[a4paper]{geometry} 
\usepackage{enumerate} 
\usepackage[sort,authoryear,round]{natbib}

\newcommand{\proofbegin}{\paragraph{Proof.}}
\newcommand{\proofend}{$\blacksquare$\bigskip}
\newtheorem{theorem}{Theorem}[section]
\newtheorem{proposition}[theorem]{Proposition}

\newtheorem{lemma}[theorem]{Lemma}

\newcommand{\eps}{\varepsilon}

\begin{document}

\title{Assessing the Distribution Consistency of Sequential Data} 

\author{Mahendra~Mariadassou and Avner~Bar-Hen}

\maketitle

\tableofcontents

\section{Introduction}
\label{Sec:introduction}

Let ${\bf X} =X_1,\ldots, X_n$ be independent observations from a repeated
experiment, and with common distribution function $F$. Let $F_n$ be the
empirical distribution and $S(X_1,\dots,X_n) = S(F_n)$ be a statistic of
the observations. The precision of $S(F_n)$ is a strictly decreasing
function of $n$ and the sample size is thus a crucial issue.

It is often possible to increase the sample size by acquiring additional
observations ${\bf X'}=X_{n+1},\dots,X_{n+k}$. This is done at additional
cost and time, for example by increasing the cohorts in clinical trials or
sequencing additional genes in molecular biology. In a parametric framework
where $F$ belongs to some family $(F_\theta)_{\theta \in \Theta}$,
$S(X_1,\dots,X_n)$ would typically be an estimator of $\theta$ satisfying
$S(F_\theta)=\theta$ and the precision, often of order $n^{-1/2}$, should
decrease by using ${\bf X'}$. However the truth is often more complex. The
use of additional observations raises at least two issues, which are
addressed in this paper. The first one is the relevance of additional
observations to the inference problem. If the additional observations ${\bf
  X'}$ do not share the distribution function $F$ with ${\bf X}$, it is
certainly unwise to expect better precision when using them in the
inference. We therefore need to assess whether ${\bf X'}$ is distributed
consistently with $F$. Focusing on the average modification induced by
extending the sample to ${\bf X'}$, we provide in Section~\ref{Sec:EdgeExp}
an approximation to the law of this modification, under the consistency
hypothesis. This approximation is then fed in Section~\ref{Sec:Test} to a
test procedure and used to control the type I error. The second issue is
the relevance of acquiring the data. If the common distribution $F'$ of
observations in ${\bf X'}$ is close to $F$, one additional observation only
is likely not to be enough to detect the difference between $F$ and
$F'$. Indeed $k$ needs to be larger than some function of $n$ for the test
to be powerful. In test language, for given $F'$ and $F$, it is similar to
finding the size sample needed to achieve a power exceeding some
threshold. This issue can be solved using results of
Section~\ref{Sec:EdgeExp} and is addressed in Section~\ref{Sec:Test}.


These two issues arise in a slightly different form in sequential tests of
hypotheses and sequential change point detection. When collecting new
observations is lengthy and costly, waiting for completion of a sample of
size $n$ before performing the analyses is not a option. In such an
instance, it is desirable to use any new observation as soon as it becomes
available. Wald's Sequential Probability Ratio Test (SRPT), introduced by
his seminal paper \citep{Wald1945} and tightly connected to the classical
Neyman-Pearson test for fixed sample size, does just this. Sequential tests
stop sampling as soon as a positive result is detected and can thus be
superior to classical tests by providing results faster than classical
tests, as the success story of the Beta-Blocker Heart Attack Trial (BHAT)
proved in 1981 when it ended 8 months earlier than scheduled with positive
results \citep{Group1981}.

But, although modifications exists to account for account for composite
hypothesis \citep{Brodsky2005}, sequential tests usually test $H_0: F =
F_0$ against $H_1: F=F_1$, \emph{i.e.} observations are either all
distributed according to $F_0$ or all distributed according to $F_1$, which
is different of our main concern, since new data can have a different
distribution function than the previous ones. Sequential change point
detection is closer in essence to our needs, although it does not perfectly
fits our need either.




Sequential change point detection is heavily used in statistical quality
control. It is used to answer three questions: has a production process ran
out of control, when did it ran out of control and what is the magnitude of
the change ? Assume that the observations are distributed according to
$F_0$ under the state of control and according to $F_1$ under the other
state.
Noting $T$ the point in time at which the jump is detected and $\nu$ the
point at which it occurs, most of the change point detection literature is
interested in minimizing $E[(T-\nu)^+]$, the average number of additional
observations needed to detect the change. This is very close to our
concern: new observations not being consistent with the previous ones is
equivalent to a process running out of control at time $n$.
The CUSUM (cumulative sum) charts use the current observation to detect
significant departures of the process from the state of control
\cite{Page1954}. \citet{Lai1995} showed that a moving average scheme
consisting of only a finite size observation window around the current
observation is asymptotically as efficient as the CUSUM if the window size
grows suitably fast to infinity. \citet{Brodsky2000} generalize this result
to a larger class of schemes. But all these methods are likelihood-based
and assume $F_0$ and $F_1$ are simple enough for log-likelihood ratio to be
easily computed. \citet{Benveniste1987} use weak convergence theory to
extend CUSUM to non-likelihood-based procedures. Their asymptotic local
approach use convergence of the rescaled sums of detection statistics to a
gaussian process. \citet{Lai1999} use another approach based on moderate
deviations to extend a Generalized Likelihood Ratio (GLR) to
non-likelihood-based detection statistics. 
We present in this paper an original non-likelihood based method to check
the consistency of a new batch of observations with previous ones. Our
method requires very little assumption about $F_0$ and $F_1$ and builds
upon a simple and intuitive idea: under the hypothesis of consistency, the
precision gain obtained when adding $k$ observations to the sample can
roughly be estimated by the precision loss induced by removing $k$
observations from the sample.

Our work is motivated by the study of DNA sequences. Organisms genomes are
sequenced gene by gene: when new genes become of interest for the
community, they are simultaneously sequenced in several organisms. Waiting
for all genes from all species to be sequenced before proceeding to an
analysis is of course not an option. The current standard is to use as many
genes as available: concatenating several genes into one supergene
increases the sample size -- here the gene length -- and implies a more
accurate analysis. Such concatenation implicitly assumes that every new
gene has the same evolutionary history as the others.  Unfortunately, there
is no certainty about that. It is well known that many mechanisms --
recombination, selective sweep, purifying or positive selection among
others \citep{Balding2007}-- lead different genes to have different
histories. When a new gene becomes available, it should thus be tested for
consistency before being included in the sample. If there is suspicion or
exterior information that the new gene do not share a common history with
the previous ones, the focus is on the minimum gene length necessary to
confidently assess the difference, as in the optimization of the change
point detection.

The issue of change point detection is hardly new but unlike most methods
available in the sequential tests literature the alternative hypothesis is
not well specified: a gene can be affected by a number of evolutionary
event and thus have a number of evolutionary histories. Specifying one, or
even a finite set, of those histories in $H_1$ is hardly better than an
educated guess. The main focus is thus on rejecting $H_0$, close in
philosophy to the Repeated Significance Test (RST)
\citep{Armitage1969,Pocock1977,O'Brien1979}. This particular issue of
assessing consistency when the alternative is not well specified can also
be found in the online learning literature and is there referred to as
concept drift \citep{Domingos2000}.

The article is organized as follows: Section~\ref{Sec:Def} introduces the
key concepts and provides intuition about the kind of results we
expect. Section~\ref{Sec:EdgeExp} present our main results, derived from
Edgeworth expansions, and discuss their strong and weak
points. Section~\ref{Sec:Test} builds upon the results of
Section~\ref{Sec:EdgeExp} to present a test of consistency of a new set of
data with previous ones. Proofs are postponed to Section~\ref{Sec:Proofs}.

\section{Definitions and Notations}
\label{Sec:Def}

\subsection{Definition of $\Delta_{n,+k}$ and $\Delta_{n,-k}$}
Let $(X_1,\ldots,X_n,\ldots)$ be a sequence of i.i.d random variables whose
common distribution function is $F_0$. Consider the sample mean for the
first $n$ terms:
\[Y_n = \frac{1}{n}\sum_{i=1}^{n} X_i\] and define:
\begin{eqnarray*}
  \Delta_{n,+k} & = & Y_{n+k}-Y_n, \\
  \Delta_{n,-k} & = & Y_{n-k}-Y_{n}.
\end{eqnarray*}

Since $\Delta_{n,+k}$ is invariant by translation of the $X_i$s, we assume
without loss of generality that the $X_k$ are centered ($E[X_1]=\mu=0$) and
furthermore note:
\[
E[X_1^2] = \sigma^2 \quad E[X_1^3]=\kappa \quad E[|X_1|^3]=\beta_3 < \infty
\]

A alternative definition of $\Delta_{n,+k}$ is
\begin{equation} 
\label{Eq:Delta}
\Delta_{n,+k}
= \frac{1}{n+k}\sum_{j=1}^k X_{n+j} - \frac{k}{n(n+k)}\sum_{j=1}^n X_{j}.
\end{equation}
$\Delta_{n,+k}$ (resp. $\Delta_{n,-k}$) is centered with distribution
function $F_+$ (resp. $F_-$) and variance
$\sigma^2_{n,+k}$ (resp. $\sigma^2_{n,-k}$) where
\begin{equation*}
\sigma^2_{n,+k}  =\frac{k\sigma^2}{n(n+k)} \quad \text{ and } \quad 
\sigma^2_{n,-k}  =\frac{k\sigma^2}{n(n-k)}
\end{equation*}
$\Delta_{n,+k}$ (resp. $\Delta_{n,-k}$) represent perturbations of the
sample mean induced by adding (resp. removing) $k$ units from the
sample. As one would expect, when $n$ increases perturbations to the sample
mean are the same no matter whether $k$ terms are added to or removed from
the sample. To formalize this intuition, we focus on the difference $F_+ -
F_-$. $F_+(x) - F_-(x)$ is convenient for at least two results: using
appropriate expansion techniques, we can get results about its order of
magnitude and $\sup_{x \in \mathbb{R}} |F_+(x) - F_-(x)|$, the quantity of
interest in Kolmogorov-Smirnoff test, is easy to calculate given some
expansion of $F_+(x) - F_-(x)$.

\subsection{Characteristic Function}
\label{Sec:CharFunc}

But, before proceeding to derivation of the expansion, we recall a few
properties of characteristic functions and use them to get insight into the
difference between $\Delta_{n,+k}$ and $\Delta_{n,-k}$.

Let $X$ be a real valued random variable with distribution function
$F_X$. Let $f_X$ be the characteristic function of $X$ defined as $f_X(t) =
E[e^{itX}] = \int_{-\infty}^{\infty} e^{itx}dF_X(x)$.



Hereafter and unless specified otherwise, we use the shorthands $f$ for
$f_X$, $f_+$ for $f_{\Delta_{n,+k}}$ and $f_-$ for
$f_{\Delta_{n,-k}}$. Thanks to Eq.~\eqref{Eq:Delta} and classical
properties of the characteristic function for independent random variables,
we have
\begin{equation}
  \label{Eq:f_plus}
  f_{+}(t) = f\left(\frac{t}{n+k}\right)^k f\left(\frac{-t}{n(n+k)}\right)^n.
\end{equation}

Taylor expansion around $0$ yields
\[
f_{-}(t) - f_{+}(t) \simeq \frac{kt^2}{\sigma^2n^2} \frac{k}{n}.
\]
where lower order terms have been omitted. Note that $Var(\Delta_{n,+k})
\sim Var(\Delta_{n,-k}) \sim \frac{k\sigma^2}{n^2}$. Normalizing
$\Delta_{n,+k}$ and $\Delta_{n,-k}$ so that they have asymptotic variance
$1$ and considering the difference between the characteristic function of
the normalized version yields
\begin{equation}
  \label{Eq:DL_of_f}
  f_{-}\left(\frac{nt}{\sqrt{k}\sigma}\right) - 
  f_{+}\left(\frac{nt}{\sqrt{k}\sigma}\right) \simeq
  \frac{kt^2}{n}.
\end{equation}
omitting again all lower order terms. Since the first order term in the
expansion of $f_{-}-f_{+}$ around $0$ is of order $k/n$ and although local
expansion provides is not enough to prove it, we expect from the inversion
theorem the difference $F_- - F_+$ to be of order $k/n$. However, in order
to achieve this result, two competing speeds need to be balanced: $k^{-1/2}$
and $k/n$. An intuitive justification follows.  It is clear that
\begin{equation}
  \label{Eq:DL_of_standardized_f}
  \frac{n}{\sqrt{k}\sigma}\Delta_{n,+k} = \left(1+\frac{k}{n}\right)^{-1}
  \frac{1}{\sigma\sqrt{k}}\sum_{j=1}^k (X_{n+j} - \bar{X}_n) 
  \quad 
  \text{and}
  \quad 
  \frac{n}{\sqrt{k}\sigma}\Delta_{n,-k} = \frac{1}{\sigma\sqrt{k}}\sum_{j=1}^k
  (X_{n-k+j} - \bar{X}_{n-k})
\end{equation}
where $\bar{X}_n$ is the empirical mean of an $n$-sample of
i.i.d. $X_j$. Since $\bar{X}_n=\mu+\mathcal{O}_P\left(n^{-1/2}\right)$, it
is clear from Eq.~\eqref{Eq:DL_of_standardized_f} that
$\frac{n}{\sqrt{k}\sigma}\Delta_{n,+k}$ can be thought of as the
standardized sum of $k$ i.i.d roughly centered random variables with
variance $1$. If $k$ goes to infinity with $n$, the speed $k^{-1/2}$ is
thus the usual speed of the central limit theorem whereas $k/n$ is the
speed of the first order difference between variance of $\Delta_{n,+k}$ and
$\Delta_{n,-k}$. Depending on the regularity of $F$ and the compared speed
of $k^{-1/2}$ and $k/n$, we can make the intuition rigorous and prove the
assertion:
\begin{equation}
  \label{Eq:ExpansionGeneral}
  F_+\left(\frac{\sqrt{k}\sigma x}{n}\right) - F_-\left(\frac{\sqrt{k}\sigma x}{n}\right) =
  \frac{x}{\sqrt{2\pi}}e^{-x^2/2}\frac{k}{n}+o\left(\frac{k}{n}\right) 
\end{equation}
uniformly in $x$. Proper formulations and proofs are provided in
Section~\ref{Sec:EdgeExp}.

Eq.~\eqref{Eq:DL_of_f} provides an asymptotic expansion of $f_+ - f_-$ in
an interval around $0$ and, although it gives some insight about the
resulting Eq.~\eqref{Eq:ExpansionGeneral}, it is not powerful enough to
derive it properly. We therefore resort to Edgeworth expansion, with an
Edgeworth series acting as a middleman between $f_+$ and $f_-$. This is the
aim of Section~\ref{Sec:EdgeExp}.

\section{Edgeworth Expansion}
\label{Sec:EdgeExp}

Edgeworth series provide an approximation of a probability distribution in
terms of its cumulants and are an improvement to the central limit
theorem. The nice property of Edgeworth expansions is that they are true
asymptotic expansions. We can thus control the error between a
probability distribution and its Edgeworth expansion. The literature about Edgeworth
expansion is quite abundant and full of powerful results. However most, if
not all, of these results rely heavily on $f$ satisfying the so-called
\emph{Cramér's Condition}:
\begin{equation}
  \label{Eq:Cramer}
  \limsup_{|t|\rightarrow \infty} |f(t)| < 1 
\end{equation}
Cramér's condition is equivalent to $F$ having an absolutely continuous
component \citep{Hall1984} but we take a special interest in non-lattice
completely discontinuous $F$ (\emph{i.e.} discrete $X$) for which
condition~\eqref{Eq:Cramer} is not satisfied. We deal with distribution
functions satisfying Cramér's condition in Section~\ref{Sec:Cramer} before
turning to non-lattice discrete distribution functions in
Section~\ref{Sec:Discrete}. Proofs are postponed in Section~\ref{Sec:Proofs}.

\subsection{With Cramér's Condition}
\label{Sec:Cramer}
 
The main result of this section is the following:
\begin{theorem} \label{Theo:ExpansionCramer} Let $(X_i)$ be a sequence of
  i.i.d. real valued random variables with distribution function
  $F$. Suppose that Cramér's condition holds, \emph{i.e.} that
  $\limsup_{|t|\rightarrow \infty} |f(t)| < 1$. Suppose furthermore that
  there exists an integer $m \geq 1$ such $E[|X|^{m+2}] < \infty$ and consider
  $\alpha \in \left(\frac{2}{m+2},1\right)$. If $k \sim n^{\alpha}$ then:
  \begin{equation}
    F_+ \left(\frac{\sqrt{k}\sigma x}{n} \right) - F_-
    \left(\frac{\sqrt{k}\sigma x}{n}\right) =
    \frac{xe^{-x^2/2}}{\sqrt{2\pi}}\frac{k}{n}+o\left(\frac{k}{n}\right)
  \end{equation}
  uniformly in $x$.
\end{theorem}

If $E[|X|^m] < \infty$ for all $m$, as is the case for gaussian random
variables, $\alpha$ can take any value in $(0,1)$. The only missing case is
$k = o(n^\eps)$ for all $\eps > 0$. In particular and unlike gaussian
variables, as will be shown in Prop.~\ref{Prop:Normal}, $k$ can not be fixed or grow only logarithmically with $n$.

\subsection{Without Cramér's Condition}
\label{Sec:Discrete}
 
The main result of this section is the following:
\begin{theorem} \label{Theo:ExpansionDiscrete} Let $(X_i)$ be a sequence of
  i.i.d. real valued random variables with distribution function
  $F$. Suppose that $X$ is a non lattice, discrete random variable. Suppose
  furthermore that $\beta_3 = E[|X|^3] < \infty$ and consider $\alpha \in
  \left(\frac{2}{3},1\right)$. If $k \sim n^{\alpha}$ then:
  \begin{equation}
    F_+ \left(\frac{\sqrt{k}\sigma x}{n} \right) - F_-
    \left(\frac{\sqrt{k}\sigma x}{n}\right) =
    \frac{xe^{-x^2/2}}{\sqrt{2\pi}}\frac{k}{n}+o\left(\frac{k}{n}\right)
  \end{equation}
  uniformly in $x$.
\end{theorem}

The fundamental difference between Theorems~\ref{Theo:ExpansionDiscrete}
and \ref{Theo:ExpansionCramer} lies in the range of value $\alpha$ can
take. When the distribution function $F$ of $X$ has some absolutely
continuous component, $k$ is allowed, upon moment conditions, to grow
slowly compared to $n$. When the distribution function is completely
discrete, the third order moment is enough to achieve the expansion. Higher
order moments, even if they do exist, are not sufficient to expand the
range of value $\alpha$ can take and are thus not required.

\subsection{New Generating Process}
\label{Sec:NewProcess}

The main result of this section is the following: 
\begin{theorem} \label{Theo:ExpansionTest} Let $X_i$ (resp. $Y_i$) be a
  sequence of i.i.d. real valued random variables with distribution
  function $F_0$ (resp. $F_1$). Suppose that $X$ (resp. $Y$) has finite
  expectation $\mu_0$ (resp. $\mu_1$) and variance $\sigma_0^2$
  (resp. $\sigma_1^2$). Suppose furthermore that $\beta_3 = E[|Y|^3] <
  \infty$ and consider $\alpha \in (0,1)$. If $k \sim n^\alpha$, then:
  \begin{equation}\label{Eq:ExpansionTest}
    F_+ \left(\frac{\sqrt{k}\sigma_1 x}{n} \right) = \Phi\left(x  -
      \frac{n}{n+k} \frac{\sqrt{k}(\mu_1 - \mu_0)}{\sigma_1}\right) +
    \mathcal{O}(n^{-\beta})
  \end{equation}
  uniformly in $x$, where $\beta = \min{(\frac{\alpha}{2},1-\alpha)}$. If
  $x$ is restricted to a bounded range and $\mu_1 \neq \mu_0$, the
  correcting term $n/(n+k)$ is unnecessary and Eq.~\eqref{Eq:ExpansionTest}
  simplifies to
  \begin{equation}
    F_+\left(\frac{\sqrt{k}\sigma_1 x}{n} \right) = \Phi\left(x  -
      \frac{\sqrt{k}(\mu_1 - \mu_0)}{\sigma_1}\right) +
    \mathcal{O}(n^{-\beta}).
  \end{equation}

\end{theorem}

Theorem~\ref{Theo:ExpansionTest} requires a third order condition on the
new generating process $Y$ to ensure that the remaining term is of order
$\mathcal{O}(k^{-1/2})$. Neglecting second order terms,
$\frac{n\Delta_{n,+k}}{\sqrt{k}\sigma_1}$ behaves like a gaussian variable
with mean $\sqrt{k}\frac{\mu_1 - \mu_0}{\sigma_1}$ and variance $1$. As we
could expect, the mean diverges faster if $\mu_0$ and $\mu_1$ are well
separated when compared to the scale $\sigma_1$.

\subsection{About Discrete Distributions}

Our motivating example of DNA analysis is intimately linked to discrete
state space. When comparing the same gene among a set of $s$ organisms,
each nucleotide in a species is associated to its homologous in the
remaining species. An observation consists of a $s$-uple of nucleotides,
. Each nucleotide can take value in the set $\{A,C,G,T\}$ and thus the
$s$-uples take value in $\{A,C,G,T\}^s$. The statistic of interest is the
likelihood of an observation under a given model. The observations are
intrinsically discrete and so is the likelihood of an observation under a
given model. To turn these likelihoods to continuous variables and allow
for the use of Theorem~\ref{Theo:ExpansionCramer} instead of the less
powerful Theorem~\ref{Theo:ExpansionDiscrete}, we must resort to the
trick exposed hereafter.

Formally, consider a discrete space $A = (a_i)_{i=1,\dots,N}$ and a
probability measure $\mathbf{\theta}=(\theta_1,\ldots,\theta_N)$ on $A$. In
DNA analysis, $A = \{A,C,G,T\}^s$ and $\mathbf{\theta}$ is a model
assigning a probability to each $a \in A$. Assume $\theta_i > 0$ for all
$i$ and let $(Z_i)_{i \in \mathbb{N}}$ be a sequence of i.i.d. random
variables such that $P(Z = a_j) = \theta_j$ for $j=1,\dots,N$. We take a special
interest in $(X_i)_{i \in \mathbb{N}}$ defined as
\[
X_i = \log{P(\{Z_i\})} = \sum_{j=1}^N \log{P(Z_i=a_j)}\mathbbm{1}_{\{Z_i =
  a_j\}}
\]
$(X_i)$ is easily an i.i.d sequence of discrete random variables such that
$P(X = \log(\theta_j)) = \theta_j$.  In this case, we can prove thanks to
Theorem~\ref{Theo:ExpansionDiscrete} that $\sup_{\mathbb{R}} |F_+ - F_-| =
\frac{1}{\sqrt{2\pi e}}\frac{k}{n} + o\left(\frac{k}{n}\right)$ but only if
$k \sim n^{\alpha}$ with $\alpha \in (2/3,1)$. We don't have access to
lower values of $\alpha$.

Suppose now that $\mathbf{\theta}$ is not the same for all $Z_i$ but rather
that each $Z_i$ is drawn from $A$ according to a specific
$\mathbf{\alpha}^{(i)}=(\alpha^{(i)}_1,\ldots,\alpha^{(i)}_N)$ and
furthermore that $\mathbf{\alpha^{(i)}}$ is an i.i.d sequence from a
Dirichlet distribution $\text{Dir}(\lambda\mathbf{\theta})$ that has
density:
\[
f(v_1,\ldots,v_{N-1}) = \frac{\prod_{i=1}^N\Gamma(\lambda
  \theta_i)}{\Gamma(\lambda\sum_{i=1}^N \theta_i)} \prod_{i=1}^{N-1}
v_i^{\lambda \theta_i -1}
\]
for all $v_1,\ldots,v_{N-1} > 0$ such that $\sum_{i=1}^{N-1} v_i < 1$ and $
V_N = 1 - \sum_{i=1}^{N-1} V_i$. Intuitively, $(V_1,\ldots,V_N)$ is a
vector of the $N$ dimensional unit simplex with mean $\mathbf{\theta}$ and
variance inversely proportional to $\lambda$: the marginal distribution of
$V_i$ has mean $\theta_i$ and variance
$\frac{\theta_i(1-\theta_i)}{\lambda+1}$. Using $Dir(\lambda
\mathbf{\theta})$ instead of $\mathbf{\theta}$ can be seen as a
regularization of the previous case, with $\mathbf{\theta}$ being the
limiting case of $Dir(\lambda\mathbf{\theta})$ when $\lambda$ goes to
infinity.

It is then easily seen that the $X_i$ are i.i.d random variables taking
value in $\mathbb{R}_-$ and absolutely continuous with respect to the
Lebesgue-measure. A bit of algebra gives for all $m$
\begin{eqnarray*}
  E[|X|^m] = E\left[\sum_{i=1}^N |log^m P(Z = a_i)| \mathbbm{1}_{Z = a_i}
  \right] & = & \sum_{i=1}^N \int_{0}^1 |\log(\alpha_i)|^m \alpha_i
  p(\alpha_i|\mathbf{\theta})d\alpha_i \\ 
  & = & \sum_{i=1}^N \frac{\Gamma(\lambda
    \theta_i)\Gamma(\lambda(1-\theta_i))}{\Gamma(\lambda)}  \int_{0}^1
  |log^m(x)|x^{\lambda \theta_i}(1-x)^{\lambda(1-\theta_i)-1}
  dx \\ 
  & < & \infty 
\end{eqnarray*}
In this case of particular interest, Theorem~\ref{Theo:ExpansionCramer}
applies for any value of $\alpha$ in $(0,1)$ as $m$ can be taken arbitrary
large.

\section{Application to Test}
\label{Sec:Test}

Theorems~\ref{Theo:ExpansionCramer} and \ref{Theo:ExpansionDiscrete} are
useful for detecting changes in the generating process of new observations.

We want to test whether the new batch of observations is generated by the
same process as the previous observations. Formally, given two
probability distributions $F_0$ and $F_1$, and a sequence of independent
random variables $(X_i)$ with associated distribution function $F_{X_i}$,
we want to test $H_0$: ``$F_{X_i} = F_0$ for $i=1,\dots,n+k$'' against
$H_1$: ``$F_{X_i} = F_0$ for $i \leq n$ and $F_{X_i} = F_1$ otherwise''.

In our problem, the statistic of interest is the sample mean, calculated
either on all $n+k$ observations ($Y_{n+k}$)or only the previous $n$
observations ($Y_n$). We shall therefore assume that $F_0$ and $F_1$ have
different means $\mu_0$ and $\mu_1$. $\Delta_{n,+k} = Y_{n+k} - Y_n$
represents the influence of the batch of $k$ new observations on the mean,
\emph{i.e} the translation of the sample mean induced by adding the batch
of new observation to the calculation. The use of the term ``influence'' is
not coincidental: $\Delta_{n,+k}$ is strongly connected to influence
functions \citep{Hampel1974,Huber2004}. When the quantity to estimate is
the mean $\mu$ of a distribution and $k=1$, $n\Delta_{n,+1}$ is indeed
exactly the empirical influence value of observation $X_{n+1}$ on the
estimator $Y_n = \frac{1}{n}\sum_{i=1}^n X_i$ of $\mu$, \emph{i.e.} the influence of an
infinitesimal perturbation on $\hat{\mu}$ along the direction
$\delta_{X_i}$, the unit mass at point $X_i$.

Large positive or negative influence values point up the corresponding
observations as potentials outliers whereas small to moderate influence
values support consistency of the data. Up to a rescaling, $\Delta_{n,+k}$
can be understood as an extension of influence functions to a batch of
observations instead of a single one.

\subsection{Distribution of $\Delta_{n,+k}$ under $H_0$}
\label{Sec:H0}

Let $k \in \{n^{\beta_1},n^{\beta_2}\}$ with $\beta_1$ and $\beta_2$ to be
specified later. Under $H_0$, $F_{X_i} = F_0$ for $i=1,\dots,n+k$ and it
comes from Theorems~\ref{Theo:ExpansionCramer} for continuous and
\ref{Theo:ExpansionDiscrete} for discrete distributions that
$\Delta_{n,+k}$ and $\Delta_{n,-k}$ have the same distribution function, up
to a correcting term of order $k/n$. For discrete distributions,
$(\beta_1,\beta_2) = (2/3+\eps, 1-\eps)$ where $\eps$ is an arbitrary small
positive value. For continuous distributions $(\beta_1,\beta_2) =
(\frac{2}{m+2}+\eps, 1-\eps)$ where $\eps$ is again an arbitrary small
positive value and $m$ is the highest order moment of $F_0$.

The alternative definition Eq.~\eqref{Eq:Delta} of $\Delta_{n,-k}$ gives
different weights to $(X_1,\dots,X_{n-k})$ and
$(X_{n-k+1},\dots,X_n)$. Under $H_0$, the first $n$ observations are
identically distributed and exchangeable. Exchangeability implies that the
order of $(X_1,\dots,X_n)$ does not matter. Since their order does not
matter, $(X_{n-k+1},\dots,X_n)$ can be replaced by any other subset of
$(X_1,\dots,X_n)$ of size $k$. In particular, the distribution of
$\Delta_{n,-k}$ can be approximated by repeatedly selecting $k$ terms from
$(X_1,\dots,X_{n})$ and substituting them to $(X_{n-k+1},\dots,X_n)$.

When the distribution $F_0$ of the $X_i$ under $H_0$ is not a simple
parametric function or involves a large number of parameters, the exact
distribution function of $\Delta_{n,+k}$ is unachievable. Even an Edgeworth
expansion \emph{à la} Prop.~\ref{Prop:f_something} requires the estimation
of many cumulants. By contrast a good numerical approximation of $F_-$ is
available thanks to the previous remark and we can substitute it to
$F_+$. Adding the correcting term of order $k/n$ only requires the
estimation of the standard deviation $\sigma$ of $F_0$. And one may notice
that since there are $n+k$ observations with $n$ larger than $k$, the
estimation of $\sigma$ is significantly more accurate than the
approximation of $F_-$ by its empirical version.

Wrapping up the preceding remarks, the distribution $F_+$ of
$\Delta_{n,-k}$ can approximated in the following way:
\begin{enumerate}[(i)]
\item Compute the mean $Y_n$ of the $n$ observations;
\item Select at random without replacement $k$ observations among the $n$;
\item Compute the mean $Y^\star_{n-k}$ of the remaining $n-k$ observations;
\item Record the difference $\Delta_{n,-k}^\star = Y_n - Y_{n-k}^\star$;
\item Repeat (ii) to (iv) a large number $(N)$ of times.
\end{enumerate}
The distribution $F_+$ of $\Delta_{n,+k}$ is then well approximated by the
distribution of $\Delta_{n,-k}^\star$, corrected by the term of order $k/n$
(see \citet{Hall1984} for more detailed results). The approximation of
$F_+$ can then be used to construct a critical region for rejecting $H_0$
based on the $\Delta_{n,+k}$.

\subsection{Distribution of $\Delta_{n,+k}$ under $H_1$}
\label{Sec:H1}

Under $H_1$, noting $\sigma^2_1$ the variance of the distribution $F_1$ and
assuming $\mu_0 \neq \mu_1$, Theorem~\ref{Theo:ExpansionTest} implies
\[
F_+\left(\frac{\sqrt{k}\sigma_1 t}{n}\right) = \Phi\left(t -
  \sqrt{k}\frac{\mu_1 - \mu_0}{\sigma_1}\right) + \mathcal{O}(k^{-1/2}) +
\mathcal{O}\left(\frac{k}{n}\right)
\]
where $\Phi$ is the standard normal distribution. The distribution of
$\Delta_{n,+k}$ under $H_1$ is approximately gaussian with mean
$\sqrt{k}\frac{\mu_1 - \mu_0}{\sigma_1}$ diverging to $\infty$ with
$k$. Difference between $F_+$ and $F_-$ is of order $\mathcal{O}(1)$ and
terms correcting for the lack of gaussianity of the observations are
negligible in front of the main term. Given the boundary of the rejection
zone calculated in section~\ref{Sec:H0}, the approximate power of the test
can then easily be computed. 




\subsection{Discussion of the results}

\paragraph{About the remainder term:} Theorems~\ref{Theo:ExpansionCramer}
and \ref{Theo:ExpansionTest} are derived for very general distribution
functions: they hold under mere moment conditions. When the distribution at
hand is better specified, more accurate results can reasonably be
expected. But in the absence of any further assumptions, the remainder of
order $o(k/n)$ is possibly the best we can achieve.

For example, if the distribution function is skewed, tedious calculations
show that the remainder is at least of order $\mathcal{O}(\sqrt{k}/n)$. And
we can get closer to $k/n$ by mimicking discrete lattice
distributions. Lattice distributions are off-limits but can be seen as the
limiting case of non-lattice discrete distributions: a discrete non-lattice
distribution with jumps of size $1/2 - \eps$ at points $\pm 1$ and size
$\eps$ at points $\pm \sqrt{2}$ is very close to a lattice distribution
with jumps of size $1/2$ at points $\pm 1$ for small enough $\eps$. For the
limiting case of $F_0$ being such a lattice distribution, and for odd $k$
such that neither $n/k$ nor $(n-k)/k$ are integer, $F_+$ has a jump of size
of asymptotic size $\sqrt{2/\pi k}$ at point $1/(n+k)$ when $F_-$ has no
jump at that point. Since $\frac{kx}{n}e^{-x^2/2}$ has no jump whatsoever
at any point, the extremum of $(F_+(\sqrt{k}\sigma x/n) -
F_-(\sqrt{k}\sigma x/n)) - kx/n e^{-\frac{x^2}{2}}$ is at least
$\sqrt{2/\pi k}$ attained for $x = \frac{n}{n+k} \frac{1}{\sqrt{k}\sigma}$
and thus of order at least $k^{-1/2}$. Since $k^{-1/2} \sim n^{-\alpha/2}$
which can be arbitrarily close to $k/n$ as $\alpha$ decreases towards
$2/3$, the $o(k/n)$ can not be improved upon in this case.

On the other hand, gaussian variables have such a nice distribution that
most calculations about $F_-$ and $F_+$ can be done exactly. Most important
of all, whatever the value of $k$, if $(X_{n+1},\dots,X_{n+k})$ is a linear
vector, then any linear combination of $X_{n+1},\dots,X_{n+k}$ is
gaussian. Going back to Eq.~\eqref{Eq:DL_of_standardized_f}, the first term
is \emph{exactly} gaussian and there is no need whatsoever for correcting
terms of order $k^{j/2}$. This is the most favorable case, for which the
remainder in Theorem~\ref{Theo:ExpansionCramer} has the smallest order of
magnitude.

Under $H_0$, if the $X_i$ have mean $\mu$ and variance $\sigma^2$, then
$\frac{n\Delta_{n,-k}}{\sqrt{k}\sigma} \sim \mathcal{N}(0,\frac{n}{n-k})$,
$\frac{n\Delta_{n,+k}}{\sqrt{k}\sigma} \sim \mathcal{N}(0,\frac{n}{n+k})$
and we can derive the following result:
\begin{proposition}
  \label{Prop:Normal}
  Let $\Delta_{n,+k}$ and $\Delta_{n,-k}$ be defined as before, then:
  \[
  F_+\left(\frac{\sqrt{k}\sigma x}{n}\right) -
  F_-\left(\frac{\sqrt{k}\sigma x}{n}\right) -
  \frac{k}{n}\frac{xe^{-\frac{x^2}{2}}}{\sqrt{2 \pi}} =
  \mathcal{O}\left(\frac{k^3}{n^3}\right)
  \]
  Uniformly in $x$.
\end{proposition}
Prop~\ref{Prop:Normal} is better than the result provided by
Theorem~\ref{Theo:ExpansionCramer}, as $\mathcal{O}(k^3/n^3)$ is smaller
than $o(k/n)$. Further algebra can even prove here that
$\mathcal{O}(k^3/n^3)$ is no greater than $1.2 k^3/n^3$, uniformly in $x$.

Under hypothesis $H_1$, we have:
\[
\frac{n\Delta_{n,+k}}{\sqrt{k}\sigma} \sim
\mathcal{N}\left(
\frac{n}{n+k}\frac{\sqrt{k}(\mu_1 - \mu_0)}{\sigma_1},
1+\frac{Ak}{n} \right)
\]
where $A \leq 1+\frac{\sigma_0^2}{\sigma_1^2}$. As expected, the result is
again slightly more accurate than would be obtained by
Theorem~\ref{Theo:ExpansionTest} alone, as the remainder is exactly,
instead of at least, of order $(k/n)^{1/2}$. In the gaussian case, we can
thus easily improve upon results from Section~\ref{Sec:EdgeExp}.

\paragraph{About Cramer's Condition:} Cramer's condition plays a crucial
role in the demonstration of Theorem~\ref{Theo:ExpansionCramer}. Without
Cramer's condition, there is no guarantee that jumps of the distribution
function $F_+$ are of order $o(k^{-1})$ and higher order moments of $F_+$
can not be used to improve the range of $k$ that can be used. Indeed, as
the binomial example emphasizes for the forbidden but limiting case of
lattice distribution, jumps can be of order $k^{-1/2}$. But for non-lattice
discrete lattice distributions, the maximum jump is at most of order
$o(k^{-1/2})$ and can be much smaller than that, for example
$o(k^{-1})$. In this case, is might be possible upon further work to
increase the range of value $\alpha$ can take in
Theorem~\ref{Theo:ExpansionDiscrete}.

\section{Proofs}
\label{Sec:Proofs}
Before we proceed to proof of Theorem~\ref{Theo:ExpansionCramer},
\ref{Theo:ExpansionDiscrete} and \ref{Theo:ExpansionTest}, we recall some
lemma concerning the expansion of $f^k(x/\sqrt{k})$.

Without loss of generality, we assume $E[X]=0$. Note $\sigma^2$ the
variance of $X$, $\alpha_j = E[X^j]$ the moment of order $j$ and $\kappa_j$
the $j$-th cumulant of $X$, defined as:
\[
\kappa_j = \left.\frac{1}{i^j}\frac{\text{d}^j}{\text{dt}^j}\ln
  E\left[e^{itX}\right]\right|_{t=0} = \frac{1}{i^j}(\ln\circ f)^{(j)}(0)
\]

\subsection{Previous Results}
\begin{lemma}[Esseen45]
  \label{Lem:Esseen1}
  Let $(X_i)$ a sequence of i.i.d. random variables and $m \geq 3$ an integer
  such that $E[|X|^m]<\infty$, then
  \[
  \left| f_X\left(\frac{t}{\sqrt{k}\sigma}\right)^k -
    e^{-\frac{t^2}{2}}\left(1 + \sum_{j=1}^{m-2} \frac{P_j(it)}{k^{j/2}}
    \right) \right| \leq \frac{\delta(k)}{k^{\frac{m-2}{2}}}
  (|t|^m+|t|^{3(m-1)})e^{-\frac{t^2}{4}} \quad \text{for} \quad |t| \leq
  \frac{\sigma\sqrt{k}}{4\beta_m^{1/m}}
  \]
  where $P_j(it) = \sum_{v=1}^j c_{jv}(it)^{2v+j}$ is a polynomial of
  degree $3j$ in $it$, the coefficient $c_{jv}$ being a polynomial in the
  cumulants $\kappa_3,\ldots,\kappa_{j-v+3}$ and $\delta(k)\rightarrow 0$.
\end{lemma}

\begin{lemma}[Esseen45]
  \label{Lem:Esseen2}
  Let $(X_i)$ a sequence of i.i.d. random variables and $2< \nu \leq 3$ a
  real number such that $\beta_\nu = E[|X|^\nu]<\infty$, then there exists
  a constant $C_\nu$ depending only on $\nu$ such that
  \[
  \left| f_X\left(\frac{t}{\sqrt{k}\sigma}\right)^k -
    e^{-\frac{t^2}{2}}\right| \leq \frac{C_\nu}{n^{\frac{\nu-2}{2}}}
  \frac{\beta_\nu}{\sigma^\nu} |t|^\nu e^{-\frac{t^2}{4}} \quad \text{for}
  \quad |t| \leq
  \frac{\sigma^{\frac{1}{\nu-2}}\sqrt{k}}{(24\beta_\nu)^{\frac{1}{\nu-2}}}
  \]
\end{lemma}

Lemma~\ref{Lem:Esseen1} and \ref{Lem:Esseen2} are proved in
{\cite{Esseen1945}} (p. 44). An alternative proof can be found in
{\cite{Cramer1937}} (p. 71 and 74).

\begin{lemma}[Esseen48]
  \label{Lem:Esseen3}
  Let $X$ be a non lattice discrete random variable, then for every $\eta >
  0$ there exists a positive function $\displaystyle \lambda(k) \underset{k
    \to \infty}{\longrightarrow} \infty$ such that:
  \[
  \int_{\eta}^{\lambda(k)} \frac{|f(t)|^k}{t} =
  o\left(\frac{1}{\sqrt{k}}\right)
  \]
\end{lemma}
The proof of Lemma~\ref{Lem:Esseen2} can be found in
{\cite{Esseen1945}} (Lemma~1, p. 49).

We recall one last theorem before proceeding to the proof.
\begin{theorem}[Essen48]
  \label{Theo:Esseen}
  Let $A,T$ and $\eps$ be arbitrary positive constants, $F(x)$ a
  non-decreasing function, $G(x)$ a real function of bounded variation on
  the real axis, $f(t)$ and $g(t)$ the corresponding Fourier-Stieltjed
  transforms such that:
  \begin{enumerate}
  \item $F(-\infty) = G(-\infty) = 0$, $F(\infty) = G(\infty)$
  \item $G'(x)$ exists everywhere and $|G'(x)| \leq A$
  \item $\int_{-T}^T \left|\frac{f(t)-g(t)}{t}\right|dt = \eps$
  \end{enumerate}
  To every number $k > 1$, there corresponds a finite positive number
  $c(k)$, only depending on $k$, such that
  \[
  |F(x) - G(x)| \leq k\frac{\eps}{2\pi} + c(k)\frac{A}{T}
  \]
\end{theorem}

The proof of Theorem~\ref{Theo:Esseen} is given in {\cite{Esseen1945}}
(Theorem 2.a, p. 32)

\subsection{New Results}
Lemma~\ref{Lem:Esseen1bis} is a generalization of Lemma~\ref{Lem:Esseen1}.

\begin{lemma}
  \label{Lem:Esseen1bis}
  Suppose that $X_i$ is a sequence of i.i.d. random variables such
  $E[|X|^m]<\infty$ for an integer $m \geq 3$, then for $|t| \leq
  \frac{\sigma\sqrt{k}}{4\beta_m^{1/m}}$:
  \begin{eqnarray*}
    \left| f_X\left(\frac{t}{\sqrt{k}\sigma}\frac{n}{n+k}\right)^k -
      e^{-\frac{t^2}{2}}\left(1+\frac{kt^2}{n}\right)\left(1 +
        \sum_{j=1}^{m-2} \frac{P_j(it)}{k^{j/2}} \right) \right| & \leq &
    \left\{\frac{\delta(k)}{k^{\frac{m-2}{2}}}+C_m\frac{\sqrt{k}}{n}\right\}
    (|t|^3 +|t|^{3(m-1)})e^{-\frac{t^2}{4}} \\
    & & \quad + C'_m \frac{k^2}{n^2}(|t|^4+|t|^{3m-2})e^{-\frac{t^2}{4 }}
  \end{eqnarray*}

  where $P_j(it) = \sum_{v=1}^j c_{jv}(it)^{2v+j}$ is a polynomial of
  degree $3j$ in $it$, the coefficient $c_{jv}$ being a polynomial in the
  cumulants $\kappa_3,\ldots,\kappa_{j-v+3}$, $\lim_{k \to \infty}\delta(k)
  = 0$ and $C_m$ and $C'_m$ are constants depending only on $m$.
\end{lemma}

\proofbegin It follows from Lemma~\ref{Lem:Esseen1} that
\[
\left| f_X\left(\frac{t}{\sqrt{k}\sigma}\frac{n}{n+k}\right)^k -
  e^{-\frac{t^2}{2}(1+\frac{k}{n})^{-2}}\left(1 + \sum_{j=1}^{m-2}
    \frac{P_j(it(1+\frac{k}{n}))}{k^{j/2}} \right) \right| \leq
\frac{\delta(k)}{k^{\frac{m-2}{2}}}
(|t|^m+|t|^{3(m-1)})\left(1+\frac{k}{n}\right)^{-3(m-1)}e^{-\frac{t^2(1+\frac{k}{n})^{-2}}{4}}
\]
We now expand $e^{-\frac{t^2}{2}(1+\frac{k}{n})^{-2}}$ in power of
$\frac{k}{n}$ and arrange the terms in a convenient order.
\[
-\frac{t^2}{2}\left\{\left(1+\frac{k}{n}\right)^{-2} - 1\right\} =
\frac{kt^2}{n} - \frac{k^2t^2}{2(n+k)^2}
\]
Furthermore $-\frac{t^2}{2} \leq
-\frac{t^2}{2}\left(1+\frac{k}{n}\right)^{-2} \leq -\frac{t^2}{4}$, where
the last inequality holds for large enough $n$. It then follows from a
Taylor expansion that
\[
\left|e^{-\frac{t^2}{2}(1+\frac{k}{n})^{-2}} -
  e^{-\frac{t^2}{2}}\left(1+\frac{kt^2}{n}\right)\right| \leq
e^{-\frac{t^2}{2}(1+\frac{k}{n})^{-2}}\left(\frac{kt^2}{n} -
  \frac{k^2t^2}{2(n+k)^2}\right)^2 \leq
e^{-\frac{t^2}{4}}\frac{k^2t^4}{n^2}.
\]
We also have, for any integer $j$
\[
(it)^j\left(1+\frac{k}{n}\right)^{-j} - (it)^j =
-(it)^j\left\{\frac{jk}{n}+\mathcal{O}\left(\frac{k^2}{n^2}\right)\right\}.
\]
And thus there exist a constant $K_j$, not depending on $n$ and $j$ such that
\[
\left| P_j(it(1+\frac{k}{n})) - P_j(it) \right|= \left| \sum_{v=1}^j
c_{jv}(it)^{2v+j}\left\{\frac{(2v+j)k}{n}+\mathcal{O}\left(\frac{k^2}{n^2}\right)\right\}
\right|
\leq 
K_j (|t|^{j+2}+|t|^{3j})\frac{k}{n}
\]
It follows that there exists a positive constant $C_m$, depending neither
on $n$ nor $k$ such that
\[
\left| 1+\sum_{j=1}^{m-2}\frac{P_j(it(1+\frac{k}{n}))}{k^{j/2}} -
  \left(1+\sum_{j=1}^{m-2}\frac{P_j(it)}{k^{j/2}} \right)\right| \leq
\sum_{j=1}^{m-2}K_j\frac{k}{n}\frac{(|t|^{j+2}+|t|^{3j})}{k^{j/2}}  \leq
\frac{C_m}{3}\frac{\sqrt{k}}{n} (|t|^3+|t|^{3(m-2)}).
\]
Finally $e^{-\frac{t^2}{2}}\left(1+\frac{kt^2}{n}\right) \leq
3e^{-\frac{t^2}{4}}$ and there exists a constant $C'_m$ such that $\left|
  1+\sum_{j=1}^{m-2}\frac{P_j(it)}{k^{j/2}} \right| \leq
C'_m(1+|t|^{3(m-2)})$.  For any four reals $A,B,a,b$, $|AB - ab| \leq
|A(B-b)| + |b(A-a)|$. Using $A =
e^{-\frac{t^2}{2}}\left(1+\frac{kt^2}{n}\right)$, $a=
e^{-\frac{t^2}{2}(1+\frac{k}{n})^{-2}}$, $B =
1+\sum_{j=1}^{m-2}\frac{P_j(it(1+\frac{k}{n}))}{k^{j/2}}$ and
$b=1+\sum_{j=1}^{m-2}\frac{P_j(it)}{k^{j/2}}$ we obtain:
\begin{eqnarray*}
  \left| e^{-\frac{t^2}{2}(1+\frac{k}{n})^{-2}}\left(1 + \sum_{j=1}^{m-2}
      \frac{P_j(it(1+\frac{k}{n}))}{k^{j/2}} \right) -
    e^{-\frac{t^2}{2}}\left(1+\frac{kt^2}{n}\right)
    \left(1+\sum_{j=1}^{m-2}\frac{P_j(it)}{k^{j/2}}\right) 
  \right| & \leq &
  C_m\frac{\sqrt{k}}{n}(|t|^3+|t|^{3(m-2)})e^{-\frac{t^2}{4}} \\
  & & \quad + C'_m\frac{k^2}{n^2}(|t|^4+|t|^{3m-2})e^{-\frac{t^2}{4}} \\ 
\end{eqnarray*}
From which the result immediately follows. \proofend



\begin{lemma}
  \label{Lem:Proof1.3}
  With the notations previously defined and under the conditions of
  Theorem~\ref{Theo:ExpansionCramer}
  \begin{eqnarray*}
    \left| f_X\left(\frac{-\sqrt{k}t}{(n-k)\sigma}\right)^{n-k} - \left(1 -
        \frac{k}{n}\frac{t^2}{2}\right) \right| & \leq &
    K_-(t^2+t^4)\frac{k^2}{n^2} \\
    \left| f_X\left(\frac{-\sqrt{k}t}{(n+k)\sigma}\right)^{n} - \left(1 -
        \frac{k}{n}\frac{t^2}{2}\right) \right| & \leq &
    K_+(t^2+t^4)\frac{k^2}{n^2} \\
  \end{eqnarray*}
  uniformly for $|t| \leq \frac{\sigma\sqrt{k}}{4\beta_m^{1/m}}$, where
  $K_+$ and $K_-$ are constants not depending on $n,k$ or $X$.
\end{lemma}

\proofbegin Since the two inequalities are proved in the same way, we prove
only the first one. It is readily observed that $\beta_m^{1/m}$ increases
with $m$, thus $\beta_3 \leq \beta_m^{3/m}$. It follows by taking $\nu=3$
in Lemma~\ref{Lem:Esseen2} that for $|t| \leq
\frac{\sigma\sqrt{k}}{4\beta_3^{1/3}}$,
\[
\left| f_X\left(\frac{t}{\sqrt{k}\sigma}\right)^k - e^{-\frac{t^2}{2}}
\right| \leq C_3 \frac{\beta_3}{\sigma^3} \frac{1}{k^{1/2}}|t|^3
e^{-\frac{t^2}{4}}
\]
A simple decomposition of the quantity to upper bound yields
\begin{eqnarray}
  \left| f_X\left(\frac{-\sqrt{k}t}{(n-k)\sigma}\right)^{n-k} - \left(1 -
      \frac{k}{n}\frac{t^2}{2}\right) \right| & \leq & \left|
    f_X\left(\frac{-\sqrt{k}t}{(n-k)\sigma}\right)^{n-k} -
    \exp\left\{-\frac{kt^2}{2(n-k)}\right\} \right| + \label{Eq:Lemma6.1}\\ 
  & & \hspace{-0.5cm} \left| \exp\left\{-\frac{kt^2}{2(n-k)}\frac{t^2}{2}\right\} -
    \exp\left\{-\frac{kt^2}{2n}\right\} \right| + \left|
    \exp\left\{-\frac{kt^2}{2n}\right\} - \left(1 -
      \frac{kt^2}{2n}\right) \right| \notag 
\end{eqnarray}
For large enough $n$, $k \leq n-k$ and thus for $|t| \leq
\frac{\sigma\sqrt{k}}{4\beta_3^{1/3}} \leq
\frac{\sigma\sqrt{n-k}}{4\beta_3^{1/3}}$, the first term of the right-hand
side of Eq.~\eqref{Eq:Lemma6.1} is upper bounded by
\[
\left| f_X\left(\frac{-\sqrt{k}t}{(n-k)\sigma}\right)^{n-k} -
  e^{-\frac{k}{n-k}\frac{t^2}{2}} \right| \leq
C_3\frac{\beta_3}{\sigma^3}\frac{k^{3/2}}{(n-k)^{2}} |t|^3
e^{-\frac{k}{n-k}\frac{t^2}{4}} \leq K_2 \frac{k^{3/2}}{n^2}|t|^3
\]
where $K_2 = C_3\beta_3/\sigma^3\sup_n\{n^2/(n-k)^2\}$.

Using the classical inequality $|e^{x+y} - e^x| \leq |y|e^x$ for $y<0$ we
bound the second term of Eq.~\ref{Eq:Lemma6.1}:
\[
\left| e^{-\frac{k}{n-k}\frac{t^2}{2}} - e^{-\frac{k}{n}\frac{t^2}{2}}
\right| \leq e^{-\frac{k}{n}\frac{t^2}{2}}\left| \frac{k}{n-k} -
  \frac{k}{n}\right| \frac{t^2}{2} \leq K_1\frac{k^2}{n^2}t^2
\]
where $K_1 = sup_n \{ n/(n-k)\}/2$.  Finally we bound the third term of
Eq.~\ref{Eq:Lemma6.1} using the inequality $|e^{-x} - (1-x)| \leq x^2/2$
for $x \geq 0$:
\[
\left| e^{-\frac{k}{n}\frac{t^2}{2}} - \left(1 -
    \frac{k}{n}\frac{t^2}{2}\right) \right| \leq
\frac{k^2}{n^2}\frac{t^4}{4}
\]
Since $k^{3/2}/n^2 = o(k^2/n^2)$, for $K_+$ large enough
\[
K_2 \frac{k^{3/2}}{n^2}|t|^3 + K_1 \frac{k^2}{n^2}t^2 +
\frac{1}{4}\frac{k^2}{n^2}t^4 \leq K_+(t^2 + t^4)\frac{k^2}{n^2}
\]
which ends the proof of the first part of the lemma. Replacing $n-k$ by
$n+k$, the same demonstration holds and yields the second inequality of the
lemma. \proofend

\begin{lemma}
  \label{Lem:Proof1.4}
  With the notations previously defined and under the conditions of
  Theorem~\ref{Theo:ExpansionCramer}
  \begin{eqnarray*}
    \left| f_-\left(\frac{nt}{\sqrt{k}\sigma}\right) -
      e^{-\frac{t^2}{2}} \left(1-\frac{kt^2}{2n}\right) 
      \left(1 + \sum_{j=1}^{m-2} \frac{P_j(it)}{k^{j/2}}\right) 
    \right| & 
    \leq & 
    \frac{\delta(k)}{k^{\frac{m-2}{2}}}
    (|t|^3+|t|^{3(m-1)})e^{-\frac{t^2}{4}} + K'_- \frac{k^2}{n^2} e^{-\frac{t^2}{4}}(|t|^2+|t|^{3m-2})
    \\
    \left| f_+\left(\frac{nt}{\sqrt{k}\sigma}\right) -
      e^{-\frac{t^2}{2}}\left(1+\frac{kt^2}{2n}\right)\left(1 + \sum_{j=1}^{m-2}
        \frac{P_j(it)}{k^{j/2}} \right) \right| & \leq &
    \left\{\frac{\delta(k)}{k^{\frac{m-2}{2}}} +
      C_m\frac{\sqrt{k}}{n}\right\}(|t|^3+|t|^{3(m-1)})e^{-\frac{t^2}{4}}
    \\
    & & \quad + K'_+ \frac{k^2}{n^2} e^{-\frac{t^2}{4}}(|t|^2+|t|^{3m-2 })
  \end{eqnarray*}
  uniformly for $|t| \leq \frac{\sigma\sqrt{k}}{4\beta_m^{1/m}}$, where
  $K'_-$ and $K'_+$ are constants not depending on $n$ and $k$.
\end{lemma}

\proofbegin For any four reals $A,B,a,b$, $|AB - ab| \leq |B(A-a)| +
|a(B-b)|$. We take $A = f_X\left(\frac{t}{\sqrt{k}\sigma}\right)^k$, $a =
e^{-\frac{t^2}{2}}\left(1 + \sum_{j=1}^{m-2} \frac{P_j(it)}{k^{j/2}}
\right)$, $B = f_X\left(\frac{-t}{\sqrt{k}(n-k)\sigma}\right)^{n-k}$, $b =
\left(1 - \frac{kt^2}{2n}\right)$. Using $|a| \leq
C_me^{-\frac{t^2}{2}}(1+|t|^{3(m-2)})$ and Lemma~\ref{Lem:Proof1.3},
\[
|a(B-b)| \leq C_m K_- (1+|t|^{3(m-2)})(t^2+t^4)e^{-\frac{t^2}{2}} \leq
K'_-\frac{k^2}{n^2}(|t|^2+|t|^{3m-2})e^{-\frac{t^2}{4}}
\]
where $\displaystyle K'_- = C_m K_-
\sup_t\left\{e^{-\frac{t^2}{4}}\frac{|t|^2+|t|^4+|t|^{3m-4}+|t|^{3m-2}}{|t|^2+|t|^{3m-2}}\right\}$.
Similarly using $|B|\leq 1$ and Lemma~\ref{Lem:Esseen1}
\[
|B(A-a)| \leq
\frac{\delta(k)}{k^{\frac{m-2}{2}}}(|t|^m+|t|^{3(m-1)})e^{-\frac{t^2}{4}}
\]
Combining these two inequalities gives the result for the first part of the
lemma. The second part is proved in the same way using
Lemma~\ref{Lem:Esseen1bis} instead of \ref{Lem:Esseen1}. \proofend

\subsection{Proof of Prop~\ref{Prop:Normal}}
%

\begin{lemma}
  \label{Lem:Normal}
  Let $\Phi_a$ (resp. $\Phi_b$) be the cumulative distribution function of
  a centered normal random variable with variance $a$
  (resp. $b$). Furthermore assume there is $\eps >0$ such that
  $a=(1+\eps)^{-1}$ and $b=(1-\eps)^{-1}$. Then, for vanishing $\eps$:
  \[
  \Phi_a(x) - \Phi_b(x) = \eps \frac{xe^{-\frac{x^2}{2}}}{\sqrt{2\pi}} +
  \mathcal{O}\left(\eps^3\right)
  \]
  uniformly in $x$.
\end{lemma}

\proofbegin Since $\Phi_{\sigma^2}(x) = \Phi(x/\sigma)$, we have $\Phi_a(x)
= \Phi(x/\sqrt{a})$. By hypothesis $a^{-1/2} = (1+\eps)^{1/2} =
1+\eps/2-\eps^2/8+\mathcal{O}\left(\eps^3\right)$. A Taylor expansion
around 
$x$ gives
\begin{eqnarray*}
  \Phi\left(\frac{x}{\sqrt{a}}\right) & = & \Phi(x) + \Phi'(x)x(a^{-1/2}-1) +
  \frac{\Phi''(x)}{2}x^2(a^{-1/2}-1)^2 + \frac{\Phi^{(3)}(c)}{6}x^3(a^{-1/2}-1)^3 \\
  & = & \Phi(x) + x\Phi'(x)\left(\frac{\eps}{2}-\frac{\eps^2}{8}+\mathcal{O}(\eps^3)\right) +
  x^2\frac{\Phi''(x)}{2}\left(\frac{\eps^2}{4}+\mathcal{\eps^3}\right) +
  x^3\Phi^{(3)}(c)\mathcal{O}(\eps^3) \\
\end{eqnarray*}
where $c$ belongs to $(x,x/\sqrt{a})$. Since $x\Phi'(x)$ and $x^2\Phi''(x)$
can each be written $P(x)e^{-\frac{x^2}{2}}$ with $P$ a polynomial of
degree lower than $4$, they are bounded on $\mathbb{R}$. The same holds for
$x^3\Phi^{(3)}(c)$ since $|x^3\Phi^{(3)}(c)| \leq \sup_{x \in \mathbb{R}}
|x^3\Phi^{(3)}(x/\sqrt{a})| \leq 1.2a^{3/2} \leq \infty$. We can therefore
rewrite
\[
\Phi\left(\frac{x}{\sqrt{a}}\right) = \Phi(x) +
x\Phi'(x)\left(\frac{\eps}{2}-\frac{\eps^2}{8}\right) +
x^2\frac{\Phi''(x)}{2}\frac{\eps^2}{4} + \mathcal{O}(\eps^3)
\]
uniformly in $x$.  The same arguments lead to
\[
\Phi\left(\frac{x}{\sqrt{b}}\right) = \Phi(x) +
x\Phi'(x)\left(\frac{-\eps}{2}-\frac{\eps^2}{8}\right) +
x^2\frac{\Phi''(x)}{2}\frac{\eps^2}{4} + \mathcal{O}(\eps^3)
\]
Combining these two equations and using $x\Phi'(x) =
\frac{xe^{-\frac{x^2}{2}}}{\sqrt{2\pi}}$ gives the results. \proofend

\emph{Proof of Prop.~\ref{Prop:Normal}:} Since, $\frac{n^2}{{k}\sigma^2}
\sigma^2_{n,+k} = \left(1+\frac{k}{n}\right)^{-1}$ and
$\frac{n^2}{{k}\sigma^2} \sigma^2_{n,-k} =
\left(1-\frac{k}{n}\right)^{-1}$, the result is a direct consequence of
Lemma~\ref{Lem:Normal} when replacing $\eps$ by $\frac{k}{n}$.  \proofend

\subsection{Proof of Theorem~\ref{Theo:ExpansionCramer}}
\begin{proposition}
  \label{Prop:f_something}
  With the notations and under the conditions of
  Theorem~\ref{Theo:ExpansionCramer}
  \begin{eqnarray*}
    F_-\left(\frac{\sqrt{k}\sigma x}{n}\right) & = & \Phi(x) -
    \frac{kx}{2\sqrt{2\pi}n}e^{-\frac{x^2}{2}} +
    \sum_{j=1}^{m-2}\frac{P_j(-D)}{k^{j/2}}\Phi(x) + o\left(\frac{k}{n}\right)
    \\
    F_+\left(\frac{\sqrt{k}\sigma x}{n}\right) & = & \Phi(x) +
    \frac{kx}{2\sqrt{2\pi}n}e^{-\frac{x^2}{2}} +
    \sum_{j=1}^{m-2}\frac{P_j(-D)}{k^{j/2}}\Phi(x)
    + o\left(\frac{k}{n}\right)
  \end{eqnarray*}
  Uniformly in $x$, where $D$ is the differential operator.
\end{proposition}

\proofbegin The two developments are obtained in the same way, we focus on
the first one. It follows from Lemma~\ref{Lem:Proof1.4} that
\begin{eqnarray}
  A & = & \int_{-\frac{\sigma\sqrt{k}}{4\beta_m^{1/m}}}^{-\frac{\sigma\sqrt{k}}{4\beta_m^{1/m}}}
  \left| \frac{f_-\left(\frac{{nt}}{\sqrt{k}\sigma}\right) -
      e^{-\frac{t^2}{2}} \left(1-\frac{kt^2}{2n}\right)
      \left(1+\sum_{j=1}^{m-2}\frac{P_j(it)}{k^{j/2}}\right)}{t}  
  \right|dt \notag \\
  & \leq & \frac{\delta(k)}{k^{\frac{m-2}{2}}}\int_{-\infty}^\infty
  (|t|^m+|t|^{3(m-1)})e^{-\frac{t^2}{4}} + 
  \frac{K'_-k^2}{n^2}\int_{-\infty}^\infty
  (|t|^2+|t|^{3m-2)})e^{-\frac{t^2}{4}} \notag\\
  & = & \mathcal{O}\left(\frac{\delta(k)}{k^{\frac{m-2}{2}}}\right) +
  \mathcal{O}\left(\frac{k^2}{n^2}\right) = 
  o\left(\frac{k}{n}\right) \label{Eq:f_minus_around0}
\end{eqnarray}
Since Cramér's condition holds, $\sup_{|t|\geq \beta_m^{-1/m}} |f_X(t)|
\leq c < 1$.  It follows then that
\begin{equation}
  \label{Eq:f_minus_expansion}
  \int_{\frac{\sigma\sqrt{k}}{4\beta_m^{1/m}}}^{\sigma k^{\frac{m}{2}}}
  \frac{\left|f_+\left(\frac{{nt}}{\sqrt{k}\sigma}\right) \right|}{t} dt
  \leq  \int_{\frac{\sigma\sqrt{k}}{4\beta_m^{1/m}}}^{\sigma k^{\frac{m}{2}}} 
  \frac{\left|f\left(\frac{{t}}{\sqrt{k}\sigma}\right) \right|^k}{t} dt 
  \leq  \int_{\frac{1}{4\beta_m^{1/m}}}^{k^{\frac{m-2}{2}}} 
  \frac{\left|f\left(t\right) \right|^k}{t} dt 
  \leq  k^{\frac{m-2}{2}}c^k 
  =  o\left(\frac{k}{n}\right)
\end{equation}
The same holds for $e^{-\frac{t^2}{2}} \left(1-\frac{kt^2}{2n}\right)
\left(1+\sum_{j=1}^{m-2}\frac{P_j(it)}{k^{j/2}}\right)$.  Finally,
combining Eq.~\eqref{Eq:f_minus_around0} and
Eq.~\eqref{Eq:f_minus_expansion} gives
\[
\int_{-k^{m/2}}^{k^{m/2}} \left|
  \frac{f_-\left(\frac{\sqrt{nt}}{\sqrt{k}\sigma}\right) -
    e^{-\frac{t^2}{2}} \left(1-\frac{kt^2}{2n}\right)
    \left(1+\sum_{j=1}^{m-2}\frac{P_j(it)}{k^{j/2}}\right)}{t} \right|dt =
o\left(\frac{k}{n}\right)
\]
Remark that $k^{-m/2} \sim n^{-\frac{m\alpha}{2}} = o(n^{-\frac{m}{m+2}}) =
o\left(\frac{k}{n}\right)$. Using Theorem~\ref{Theo:Esseen} with $T =
k^{-m/2}$, we obtain:
\[
F_-\left(\frac{\sqrt{k}\sigma x}{n}\right) =
\left(1+\frac{kX^2}{2n}\right)(-D)\Phi(x) +
\sum_{j=1}^{m-2}\frac{\left(1+\frac{kX^2}{n}\right)P_j(-D)}{k^{j/2}}\Phi(x)
+ o\left(\frac{k}{n}\right)
\]
The term $\left(1+\frac{kX^2}{2n}\right)(-D)\Phi(x)$ of the right-hand side
gives $\Phi(x)-\frac{kx}{2\sqrt{2\pi}n}e^{-\frac{x^2}{2}}$ when doing the
inverse Fourier transform. The result then follows from
$\frac{kX^2}{n}\frac{P_j(-D)}{k^{j/2}}\Phi(x) = o\left(\frac{k}{n}\right)$
uniformly in $x$. Replacing $1+\frac{kt^2}{2n}$ with $1-\frac{kt^2}{2n}$ in
the proof gives the second expansion.

\proofend

\emph{Proof of Theorem~\ref{Theo:ExpansionCramer}:} The result is a direct
consequence from Prop.~\ref{Prop:f_something}. \proofend

\subsection{Proof of Theorem~\ref{Theo:ExpansionDiscrete}}

\textbf{Remark:} Cramér's condition is essential to ensure that the
Edgeworth expansion of $f_+$ is valid up to the order $m$. If it does not
hold, then Lemma~\ref{Lem:Esseen1} and \ref{Lem:Esseen1bis} are still valid
but $\int_\omega^t \frac{|f(t)|^k}{t}$ does not decrease exponentially fast
anymore. We are limited to $T$ or order $k^{1/2}$ in
Theorem~\ref{Theo:Esseen} so that only expansions of order $1$ are
available. But order $1$ is not enough if $n$ grows too fast compared to
$k$.

\begin{proposition}
  \label{Prop:f_something2}
  With the notations and under the conditions of
  Theorem~\ref{Theo:ExpansionDiscrete}
  \begin{eqnarray*}
    F_-\left(\frac{\sqrt{k}\sigma x}{n}\right) & = & \Phi(x) -
    \frac{kx}{2\sqrt{2\pi}n}e^{-\frac{x^2}{2}} + 
    \frac{P_1(-D)}{k^{1/2}}\Phi(x) + o\left(\frac{k}{n}\right)
    \\
    F_+\left(\frac{\sqrt{k}\sigma x}{n}\right) & = & \Phi(x) +
    \frac{kx}{2\sqrt{2\pi}n}e^{-\frac{x^2}{2}} + \frac{P_1(-D)}{k^{1/2}} \Phi(x)
    + o\left(\frac{k}{n}\right)
  \end{eqnarray*}
  Uniformly in $x$.
\end{proposition}

\proofbegin As for Prop.~\ref{Prop:f_something}, the result is an
application of Theo.~\ref{Theo:Esseen}. It follows from
Lemma~\ref{Lem:Proof1.4} that
\begin{eqnarray}
  \label{Eq:f_discrete_around0}
  A & = & \int_{-\frac{\sigma\sqrt{k}}{4\beta_3^{1/3}}}^{-\frac{\sigma\sqrt{k}}{4\beta_3^{1/3}}}
  \left| \frac{f_-\left(\frac{nt}{\sqrt{k}\sigma}\right) -
      e^{-\frac{t^2}{2}} \left(1-\frac{kt^2}{2n}\right)
      \left(1+\frac{P_j(it)}{k^{1/2}}\right)}{t} \right|dt \notag \\
  & \leq & \frac{\delta(k)}{k^{1/2}}\int_{-\infty}^\infty
  (|t|^3+|t|^6)e^{-\frac{t^2}{4}} + K'_-\frac{k^2}{n^2}\int_{-\infty}^\infty
  (|t|^2+|t|^7)e^{-\frac{t^2}{4}} \notag \\
  & = & o(k^{-1/2}) + o\left(\frac{k}{n}\right)
\end{eqnarray}

Remark that, since $\alpha > 2/3$, $k^{-1/2} \sim n^{-\alpha/2} =
o(n^{\alpha-1}) = o\left(\frac{k}{n}\right)$. Since Cramér's condition does
not hold, we resort to Lem.~\ref{Lem:Esseen3} from which it follows that
\begin{equation}
  \label{Eq:f_discrete_expansion}
  \int_{\frac{\sigma \sqrt{k}}{4\beta_3^{1/3}}}^{\sigma \sqrt{k} \lambda(k)}
  \frac{\left|f_+(\frac{nt}{\sqrt{k}\sigma})\right|}{t}dt \leq 
  \int_{\frac{\sigma \sqrt{k}}{4\beta_3^{1/3}}}^{\sigma \sqrt{k} \lambda(k)}
  \frac{\left|f(\frac{nt}{\sqrt{k}\sigma})\right|^k}{t}dt \leq 
  \int_{\frac{1}{4\beta_3^{1/3}}}^{\lambda(k)} \frac{\left|f(t)\right|}{t}dt =
  o(k^{-1/2}) = o\left(\frac{k}{n}\right)
\end{equation}
And the same holds for $e^{-\frac{t^2}{2}} \left(1-\frac{kt^2}{2n}\right)
\left(1+\frac{P_1(it)}{k^{1/2}}\right)$. Combining
Eq.~\eqref{Eq:f_discrete_around0} and \eqref{Eq:f_discrete_expansion}
yields
\[
\int_{-\sigma\sqrt{k}\lambda(k)}^{\sigma\sqrt{k}\lambda(k)} \left|
  \frac{f_-\left(\frac{nt}{\sqrt{k}\sigma}\right) - e^{-\frac{t^2}{2}}
    \left(1-\frac{kt^2}{2n}\right)
    \left(1+\frac{P_j(it)}{k^{1/2}}\right)}{t} \right|dt =
o\left(\frac{k}{n}\right)
\]
Since $\lambda(k)\rightarrow \infty$ as $k$, or equivalently $n$, goes to
infinity, $\frac{1}{\sqrt{k}\lambda(k)} = o(k^{-1/2}) =
o\left(\frac{k}{n}\right)$. Theo.~\ref{Theo:Esseen} then implies:
\[
F_-\left(\frac{\sqrt{k}\sigma x}{n}\right) =
\left(1+\frac{kX^2}{2n}\right)(-D)\Phi(x) +
\frac{\left(1+\frac{kX^2}{n}\right)P_1(-D)}{k^{1/2}}\Phi(x) +
o\left(\frac{k}{n}\right)
\]
where $D$ is the differential operator. The first term of the right-hand
side gives $\Phi(x)-\frac{kx}{2\sqrt{2\pi}n}e^{-\frac{x^2}{2}}$. The result
then follows from $\frac{kX^2}{n}\frac{P_j(-D)}{k^{j/2}}\Phi(x) =
o\left(\frac{k}{n}\right)$ uniformly in $x$. Replacing $1+\frac{kt^2}{2n}$
with $1-\frac{kt^2}{2n}$ in the proof gives the second expansion. \proofend

\emph{Proof of Theorem~\ref{Theo:ExpansionDiscrete}:} The result is a
direct consequence from Prop.~\ref{Prop:f_something2}. \proofend

\subsection{Proof of Theorem~\ref{Theo:ExpansionTest}}

\textbf{Remark:} As soon as $X$ and $Y$ have different expectations,
$\Delta_{n,+k}$ is not centered anymore and the central limit theorem is
enough to get the first order expansion of its distribution function. Up to
a normalization constant, $\Delta_{n,+k}$ drifts away to $\pm \infty$,
depending on the sign of $\mu_1 - \mu_0$.

Following along the same lines as the proofs of
Theorem~\ref{Theo:ExpansionDiscrete}, we first note that
\[
f_+\left(\frac{nt}{\sqrt{k}\sigma_1}\right) 
\exp\left\{-it\frac{n}{n+k}\frac{\sqrt{k}(\mu_1
    - \mu_0)}{\sigma_1} \right\} =
f_{Y-\mu_1}\left(\frac{n}{n+k}\frac{t}{\sqrt{k}\sigma_1}\right)^k
f_{X-\mu_0} \left(\frac{\sqrt{k}t}{(n+k)\sigma_1} \right)^n
\]
Using Lemma~\ref{Lem:Proof1.3},
\[
\left| f_{X-\mu_0} \left(\frac{\sqrt{k}t}{(n+k)\sigma_1} \right)^n -
  \left(1 - \frac{\sigma_0^2kt^2}{2\sigma_1^2n}\right) \right| \leq
K(t^2+t^4)\frac{k^2}{n^2} \quad \text{for} \quad |t| \leq \frac{\sigma_0
  \sqrt{n}}{4\beta_3^{1/3}}
\]
And it comes from Lemma~\ref{Lem:Esseen1} that
\[
\left| f_{Y-\mu_1}\left(\frac{n}{n+k}\frac{t}{\sqrt{k}\sigma_1}\right)^k -
  e^{-t^2/2}\left(1 +\frac{kt^2}{n}\right)\left(1 +
    \frac{\alpha_3}{6\sigma_1^6}\frac{(it)^3}{k^{1/2}}\right)\right| \leq
\frac{\delta(k)}{k^{1/2}}(|t|^3+|t|^6)e^{-\frac{t^2}{4}} \quad \text{for}
\quad |t| \leq \frac{\sigma_1 \sqrt{k}}{4\beta_3^{1/3}}
\]
Using the trick $|AB - ab| \leq |A(B-b)| + |b(A-a)|$ with $A =
f_{X-\mu_0}\left(\frac{\sqrt{k}t}{(n+k)\sigma_1}\right)^n$, $B =
f_{Y-\mu_1}\left(\frac{n}{n+k}\frac{t}{\sqrt{k}\sigma_1}\right)^k$, $a = 1
- \frac{\sigma_0^2kt^2}{2\sigma_1^2n}$ and $b =
e^{-t^2/2}\left(1+\frac{kt^2}{n^2}\right)\left(1 +
  \frac{\alpha_3}{6\sigma_1^6}\frac{(it)^3}{k^{1/2}}\right)$, it comes from
$|A| \leq 1$ and $|b| \leq K(1+|t|^5)e^{-\frac{t^2}{2}}$ that
\[
\left| f_{X-\mu_0}\left(\frac{\sqrt{k}t}{(n+k)\sigma_1}\right)^n
  f_{Y-\mu_1}\left(\frac{n}{n+k}\frac{t}{\sqrt{k}\sigma_1}\right)^k -
  e^{-\frac{t^2}{2}}\left(1 +
    \frac{kt^2}{2n}\left(2-\frac{\sigma_0^2}{\sigma_1^2}\right)+
    \frac{\alpha_3}{\sigma_1^3}\frac{(it)^3}{k^{1/2}}\right) \right| \leq K
\left\{\frac{k^2}{n^2}+\frac{\delta(k)}{k^{1/2}}\right\}
(t^2+|t|^9)e^{-\frac{t^2}{4}}
\]
for $ |t| \leq B =
\frac{\min(\sigma_0,\sigma_1)\sqrt{k}}{4\beta_3^{1/3}}$. It then follows
that
\[
\int^{B}_{-B} \left|
  \frac{f_{X-\mu_0}\left(\frac{\sqrt{k}t}{(n+k)\sigma_1}\right)^n
    f_{Y-\mu_1}\left(\frac{n}{n+k}\frac{t}{\sqrt{k}\sigma_1}\right)^k -
    e^{-\frac{t^2}{2}}\left(1 +
      \frac{k}{n}\frac{t^2}{2}\left(2-\frac{\sigma_0^2}{\sigma_1^2}\right)+
      \frac{\alpha_3}{\sigma_1^3}\frac{(it)^3}{k^{1/2}}\right)}{t} \right|
dt
= o\left(\frac{k}{n}\right) + o(k^{-1/2}).
\]
Lemma~\ref{Lem:Esseen3} combined to Theorem~\ref{Theo:Esseen} then provides
the following result:
\begin{eqnarray*}
  F_+\left(\frac{\sqrt{k}\sigma_1 x}{n}+\frac{k}{n+k}\frac{\mu_1 - \mu_0}{\sigma_1}\right) & = & 
  P\left\{\frac{n\Delta_{n,+k}}{\sqrt{k}\sigma_1} -
    \sqrt{k}\frac{n}{n+k}\frac{\mu_1 - \mu_0}{\sigma_1} \leq x \right\} \\
  & = & \Phi(x) +
  \frac{k}{n} \left(2-\frac{\sigma_0^2}{\sigma_1^2}\right)
  \frac{xe^{-\frac{x^2}{2}}e^{-\frac{x^2}{2}}}{2\sqrt{2\pi}} +
  \frac{\alpha_3}{6\sigma_1^3}\frac{(1-x^2)e^{-\frac{x^2}{2}}}{\sqrt{2k\pi}}
  + o\left(\frac{k}{n}\right) + o(k^{-1/2}) 
  \\
  & = & \Phi(x) + \mathcal{O}(n^{-\beta})
\end{eqnarray*}
uniformly in $x$, where $\beta = \min(\frac{\alpha}{2},1-\alpha)$. In
addition, if $x$ is bounded by some $M$, we further have
\[
F_+\left(\frac{\sqrt{k}\sigma_1 x }{n}\right) = \Phi\left(x -
  \frac{\sqrt{k}(\mu_1 - \mu_0)}{\sigma_1}\right)
\]
which concludes the proof. \proofend

\bibliographystyle{abbrvnat}

\end{document}